\documentclass[12pt]{amsart}
\usepackage{amsfonts,amssymb, amscd, latexsym, graphicx, psfrag, color}

\usepackage[all]{xy}

\newtheorem{dummy}{dummy}[section]

\newtheorem{theorem}[dummy]{Theorem}

\newtheorem{conjecture}[dummy]{Conjecture}

\newtheorem{proposition}[dummy]{Proposition}
\theoremstyle{definition}

\newtheorem{remark}[dummy]{Remark}

\newcommand{\bC}{\mathbb{C}}

\newcommand{\bP}{\mathbb{P}}
\newcommand{\bQ}{\mathbb{Q}}
\newcommand{\bR}{\mathbb{R}}
\newcommand{\bZ}{\mathbb{Z}}

\newcommand{\cM}{\mathcal{M}}
\newcommand{\cN}{\mathcal{N}}

\newcommand{\cX}{\mathcal {X}}

\newcommand{\Fuk}{\mathrm{Fuk}}
\newcommand{\Sh}{\mathit{Sh}}

\newcommand{\Hom}{\mathrm{Hom}}

\newcommand{\Perf}{\mathcal{P}\mathrm{erf}}
\newcommand{\MSh}{\mathit{MSh}}

\renewcommand{\SS}{\mathit{SS}}

\renewcommand{\SS}{\mathit{SS}}

\newcommand{\p}{\mathrm{p}}

\newcommand{\X}{X}
\newcommand{\Xvee}{X^\vee}
\newcommand{\Z}{Z}
\newcommand{\Zvee}{Z^\vee}
\newcommand{\M}{M}
\newcommand{\N}{N}
\newcommand{\Mvee}{M^\vee}
\newcommand{\Nvee}{N^\vee}
\newcommand{\T}{\mathbf{T}}
\newcommand{\Tvee}{\mathbf{T}^\vee}
\newcommand{\bSi}{\mathbf{\Sigma}}

\begin{document}

\title[Polytopes and Skeleta]{Polytopes and Skeleta}

\begin{abstract}
To each simplicial reflexive polytope $\triangle \subset \bZ^{n+1}$, we attach an $n$-dimensional space $\Lambda^\infty$.  It is the Legendrian boundary of a conic Lagrangian considered in work of the authors, Fang and Liu \cite{FLTZ}, and because of this it carries a sheaf of dg categories called the ``Kashiwara-Schapira sheaf.''   We discuss some conjectures and results about the role that $\Lambda^\infty$ and the Kashiwara-Schapira sheaf should play in homological mirror symmetry.

\end{abstract}
\author{David Treumann}
\author{Eric Zaslow}
\maketitle

\section{Introduction}

Locality in the Fukaya category of an exact symplectic manifold,
and the relation of this structure to sheaf theory, 
has been studied in works of Abouzaid \cite{A}, Abouzaid-Seidel \cite{AS},
Kontsevich \cite{K-stein}, Seidel \cite{S-morse},
and, with Sibilla, the authors \cite{STZ}; and in unpublished works of
Tamarkin, Nadler, and Nadler-Tanaka.  A common thread is to identify
a singular Lagrangian ``skeleton,'' and then relate the Fukaya category to the
combinatorial structure of how this skeleton is built from smooth parts,
as described by quivers or constructible sheaves.
Typically, these skeleta are constructed from a choice of some additional 
structure placed on the exact symplectic manifold (e.g., Stein, Liouville, Weinstein, Morse),
though the Fukaya category itself should be independent of such choices.

We wish to apply these ideas to Kontsevich's homological mirror symmetry conjecture (HMS).  In this note we describe an approach to the HMS for a certain class of toric hypersurfaces considered by Batyrev and Borisov.  Locality, skeleta, and constructible sheaves appear at the large complex structure/large volume limits of these families of hypersurfaces.

Let $\triangle$ and $\triangle^\vee$ be a dual pair of reflexive polytopes.  Batyrev and
Borisov  \cite{Ba,Bo}
explain how to construct from them a pair of smooth projective Calabi-Yau varieties $Z$ and $Z^\vee$ that should be mirror to each other.  

\begin{conjecture}[Kontsevich/Batyrev-Borisov HMS for toric hypersurfaces]
\label{conj:hms1}
There is an equivalence of categories between vector bundles on $Z$
and the Fukaya category on $Z^\vee:$
$$\Perf(\Z) \cong \Fuk(\Zvee).$$
\end{conjecture}

Let $\Z^\infty$ denote the ``large complex structure limit'' (LCL) of $\Z$---it is a reducible and therefore singular complex algebraic variety with toric components.  Let $\Zvee_\infty$ denote the ``large volume limit'' (LVL) of $\Zvee$---it is an affine hypersurface in an algebraic torus.  We review the constructions of $\Z^\infty$ and $\Zvee_\infty$ in Section \ref{sec:rpams}.  We have the following variant of Conjecture \ref{conj:hms1}, proposed by Seidel \cite{S-quartic}:

\begin{conjecture}[HMS at the limit]
\label{conj:hms2}
There is an equivalence of categories
$$\Perf(\Z^\infty) \cong \Fuk_c(\Zvee_\infty)$$
where the RHS denotes the Fukaya category of compact Lagrangian branes in $\Zvee_\infty$.
\end{conjecture}

In this paper, in Sections \ref{sec:skeleta} and \ref{sec:stackyfan},  we define a topological space $\Lambda^\infty$ which ``could be'' a Lagrangian skeleton of $\Zvee_\infty$.  The space $\Lambda^\infty$ is a Legendrian submanifold of a contact boundary of the cotangent bundle of a compact torus.  The technology of Kashiwara and Schapira,
reviewed in Section \ref{sec:microlocal},
allows us to equip $\Lambda^\infty$ with a sheaf of dg categories.

We conjecture that this sheaf of categories is equivalent to the sheaf of Fukaya categories on the skeleton of $\Zvee_\infty$.  In forthcoming work with Sibilla \cite{STZ-inprogress}, we will prove that the category of global objects of this sheaf is equivalent to the category of perfect complexes on $\Z^\infty$.  If our conjecture is accurate, this would establish HMS at the LCL/LVL limits for a large class of toric hypersurfaces.

A subconjecture is that $\Lambda^\infty$ is homotopy equivalent to the
hypersurface $Z^\vee_\infty.$
As evidence for this, we study the Leray spectral sequence
of a map $\Lambda^\infty \rightarrow \partial \triangle^\vee$,
which we conjecture degenerates at $E^2$ over the rational numbers
(this is not always true over the integers). 
We have calculated the $E^2$ page in many examples by computer -- see Section \ref{sec:cc}. 
Our results agree with the Betti numbers computed by Danilov-Khovanski\v{\i} \cite{DK}.

\subsubsection*{Acknowledgments}  We thank Bohan Fang, Chiu-Chu Melissa Liu, and Nicol\`o Sibilla for many
conversations and for collaborating on related work.  We thank David Nadler and Dima Tamarkin for their insights.  The work of EZ is supported by NSF-DMS-1104779.

\section{Reflexive polytopes and mirror symmetry}
\label{sec:rpams}

Let $\M$ be a lattice, and let $\Mvee$ denote the dual lattice.  A pair of lattice polytopes $\triangle \subset \M$ and $\triangle^\vee \subset \Mvee$ are called a \emph{reflexive pair} if they are polar duals of each other:
$$
\begin{array}{rcl}
\triangle^\vee & = & \{ m^\vee \in \Mvee \mid \forall m \in \triangle \quad \langle m^\vee,m\rangle \geq -1\}\\
\triangle & = & \{ m \in \M \mid \forall m^\vee \in \triangle^\vee \quad \langle m^\vee,m \rangle \geq -1\}
\end{array}
$$

To a reflexive pair of lattice polytopes we associate a pair of toric varieties $\X$ and $\Xvee$.  Let us make the following toric notation:
\begin{itemize}
\item We set $\N = \Mvee$ and $\Nvee = \M$.  Our hope is that, for those used to ``Fulton's notation'' \cite{Fu} in toric geometry, having four names for these two lattices will \emph{reduce} confusion.
\item Write $\N_\bR = \Mvee_\bR$ and $\Nvee_\bR = \M_\bR$ for the realification of these lattices.  We will often not distinguish $\triangle \subset \M$ (resp. $\triangle^\vee \subset \Mvee$) from its convex hull in $\M_\bR$ (resp. $\Mvee_\bR$).
\item Let $\T = \N \otimes \bC^*$ be the algebraic torus whose character lattice is $\M$, and let $\Tvee = \Nvee \otimes \bC^*$ be the algebraic torus whose character lattice is $\Mvee$.
\end{itemize}

We may associate a pair of toric varieties $\X$ and $\Xvee$ to a reflexive pair of lattice polytopes. $\X$ is a compactification of $\T$ and $\Xvee$ is a compactification of $\Tvee$, and they are characterized by the following properties:  $\X$ is Gorenstein, its anticanonical bundle is ample, and the moment polytope of the anticanonical bundle is $\triangle$.

Following Batyrev and Borisov \cite{Ba,Bo}, we expect a mirror relationship between a smooth projective anticanonical hypersurface $\Z \subset \X$ and a smooth projective anticanonical hypersurface $\Zvee \subset \Xvee$.  Our purpose is to understand this mirror relationship at the large complex structure/large volume limits.

\begin{remark}
From now on, we only study:
\begin{itemize}
\item \emph{complex} geometry of $\X$, its hypersurfaces and other related spaces
\item \emph{symplectic} geometry of $\Xvee$, its hypersurfaces and other related spaces
\end{itemize}
That is, we consider just one direction of mirror symmetry at a time.
\end{remark}

Actually it is now necessary to make some corrections, as in Batyrev-Borisov's recipe the singularities of the pairs $(\X,\Z)$ and $(\Xvee,\Zvee)$ might need to be resolved.  We will restrict our attention to pairs of reflexive polytopes where a resolution process is unnecessary: from now on we assume that $\triangle$ is vertex-simplicial, and $\triangle^\vee$ is facet-simplicial.  In that case it is natural to regard $\X$ as a smooth Deligne-Mumford stack instead of a variety with Gorenstein singularities.  (We explain how to do this in detail in Sections \ref{sec:skeleta} and
\ref{sec:stackyfan}.)  The toric variety $\Xvee$ could still be quite singular, but these singularities are not visible to the hypersurface at the large volume limit.

\begin{enumerate}
\item The \emph{large complex structure limit} or LCL is the union of the toric divisors of $\X$, regarded as a Deligne-Mumford stack.  We denote it by $\Z^\infty$.  As $\triangle$ is vertex-simplicial, it has toric components with normal crossings singularities.

\item We may write any polynomial function $f:\Tvee \to \bC$ as a finite linear combination of characters.  The \emph{Newton polytope} of $f$ is the convex hull of the characters that appear; it is a lattice polytope in $\Mvee$.  The \emph{large volume limit} or LVL is any generic affine hypersurface $\Zvee_\infty$ in $\T^\vee$ whose defining equation has Newton polytope $\triangle^\vee$.
\end{enumerate}

$Z^\infty$ evidently inherits a complex structure (necessarily singular, and possibly with some stackiness) from $X$.  We endow $Z^\vee_\infty$ with a symplectic structure (necessarily exact) induced from $\sum \frac{1}{r} dr_i \wedge d\theta_i$, where $r_i,\theta_i$ are polar coordinates on $\Tvee$.

\begin{remark}
Though $\Z^\infty$ is canonically defined, $\Zvee_\infty$ depends on a choice of Laurent polynomial $f:\Tvee \to \bC$.
\end{remark}

With this notation, we can expect Conjecture \ref{conj:hms2} to hold under a simplicial assumption:

\begin{conjecture}[HMS]
\label{conj:hms3}
Suppose $\triangle$ is vertex-simplicial and $\triangle^\vee$
is facet-simplicial.  Then there is an equivalence
$$\Perf(\Z^\infty) \cong \Fuk_c(\Zvee_\infty)$$
where the right-hand side denotes the Fukaya category of compact Lagrangian branes in $\Zvee_\infty$.
\end{conjecture}

\section{Microlocal sheaf theory}
\label{sec:microlocal}

In this section we recall some notions from the microlocal sheaf theory of Kashiwara-Schapira \cite{KS}, and its relationship to Fukaya theory \cite{NZ,N}.  In Section \ref{sec:kashschasheaf}, we introduce the ``Kashiwara-Schapira sheaf'' of dg categories, which is a variant of the theory of \cite[Chapter VI]{KS}.  In Section \ref{sec:conext} we discuss how this sheaf behaves under ``conormal extension.''

If $Y$ is a real analytic manifold, we let $\Sh(Y)$ denote the dg triangulated category of sheaves on $Y$ that are constructible with respect to an arbitrary subanalytic Whitney stratification of $Y$.  To each object $F$ of $\Sh(Y)$ we can attach a conic Lagrangian subset $\SS(F) \subset T^*Y$, called the singular support of $F$.  Given a conic subset $\Lambda \subset T^*Y$ the full subcategory of $\Sh(Y)$ spanned by objects $F$ with $\SS(F) \subset \Lambda$ is triangulated, we denote it by $\Sh(Y;\Lambda)$.  

Now $\Sh(Y;\Lambda)$ is equivalent to $\Fuk(T^*Y;\Lambda),$ the Fukaya category generated by
exact Lagrangian branes $L$ with $L^\infty \subset \Lambda^\infty$ at contact infinity \cite{NZ, N}.
It is widely expected that $\Fuk(T^*Y;\Lambda)$ represents the
global sections of a \emph{sheaf} of Fukaya categories on
$T^*Y$ supported on $\Lambda$.

\subsection{The Kashiwara-Schapira sheaf of categories on conic Lagrangians and on Legendrians}
\label{sec:kashschasheaf}
Fix a manifold $Y$ and a conic Lagrangian $\Lambda$,
so that $\Fuk(T^*Y;\Lambda)$ can be computed in terms of constructible sheaves, as above.  
Microlocal sheaf theory is well-suited to studying the expected sheaf of Fukaya categories on $\Lambda$:
using constructible sheaves,
we can define a sheaf of categories with no Fukaya theory in place.

If $\Omega \subset T^*Y$ is a conic open subset, Kashiwara and Schapira define the category of ``microlocal sheaves on $\Omega$'' to be the quotient $\Sh(Y)/\Sh(Y;T^*Y - \Omega)$.  Let us denote this category by $\MSh^\p(\Omega)$.  We may regard $\Omega \mapsto \MSh^\p(\Omega)$ as a contravariant functor from the poset of conic open subsets of $T^* Y$ to the $\infty$-category of dg triangulated categories.  In other words, $\MSh^\p$ is a presheaf of dg categories on $T^*Y$ in its conic topology.

For $\Lambda \subset T^* Y$ a conic Lagrangian, we can define a smaller variant of $\MSh^\p$ we call $\MSh^\p_\Lambda$:
$$\MSh^\p_\Lambda(\Omega) = \Sh(Y;\Lambda \cup (T^* Y - \Omega))/\Sh(Y;T^*Y - \Omega)$$
This again is a presheaf of dg categories on $T^*Y$ whose value on any open set that does not meet $\Lambda$ is zero.  It follows that the sheafification of $\MSh^\p_\Lambda$ is supported on $\Lambda$, and we can regard it as a sheaf of dg categories on $\Lambda$ itself in its conic topology.  We denote this sheaf of dg categories by $\MSh_\Lambda$.

\begin{remark}
Of course $\MSh^\p_\Lambda$ and $\MSh_\Lambda$ give different categories over open subsets of $T^* Y$, but it is not difficult to show that the sheafification process does not change \emph{global} sections: $\MSh_\Lambda(T^*Y) \cong \MSh^\p_\Lambda(T^* Y) := \Sh(Y;\Lambda)$.
\end{remark}

\begin{remark}
\label{rem:leg}
If $\Lambda$ is a conic Lagrangian let us write $\Lambda^\infty$ for the associated Legendrian subset of the contact boundary of $T^* Y$.  That is, $\Lambda^\infty = (\Lambda - (\Lambda \cap Y))/\bR_{>0}$.  As $\MSh_\Lambda$ is a sheaf in the conic topology on $\Lambda$ it induces a sheaf on $\Lambda^\infty$, which we will denote by $\MSh_{\Lambda^\infty}$.
\end{remark}

We expect that the Kashiwara-Schapira sheaf $\MSh_{\Lambda^\infty}$
agrees with the anticipated Fukaya sheaf of categories,
though to our knowledge the latter object has not yet been constructed.

\subsection{Conormal extension}
\label{sec:conext}

Let $Y$ be a manifold and let
$i:Y' \hookrightarrow Y$ be a closed submanifold.  We have a short exact sequence of vector bundles on $Y'$
$$0 \to T_{Y'}^* Y \to T^* Y \vert_{Y'} \stackrel{\pi}{\to} T^* Y' \to 0$$
We may regard the middle term of this sequence as a closed subset of $T^* Y$.  If $\Lambda$ is a conic Lagrangian in $T^* Y'$, then $\pi^{-1}(\Lambda)$ is a conic Lagrangian of $T^* Y$ that we call the \emph{full conormal extension} of $\Lambda$.

The functor $i_*:\Sh(Y') \to \Sh(Y)$ carries $\Sh(Y';\Lambda) \to \Sh(Y;\pi^{-1}(\Lambda))$---in fact this latter functor is an equivalence.  It turns out that this functor induces an equivalence between the Kashiwara-Schapira sheaves on $\Lambda$ and on $\pi^{-1}(\Lambda).$

\begin{proposition}[\cite{STZ-inprogress}]
The extension-by-zero functor $i_*$ induces an equivalence
$$\pi^{-1} \MSh_\Lambda \to \MSh_{\pi^{-1}(\Lambda)}$$
\end{proposition}

A basic consequence is the following.  If $\sigma^\circ$ is a bundle of open convex cones of the conormal bundle of $Y'$ in $Y$, then $\sigma^\circ \cap \pi^{-1}(\Lambda)$ is an open subset of $\pi^{-1}(\Lambda)$, which we call the \emph{conormal extension of $\Lambda$ determined by $\sigma^\circ$}.  The Proposition and basic properties of pullback sheaves give us an identification
$$\MSh_{\pi^{-1}(\Lambda)}(\Omega) \cong \Sh(Y';\Lambda)$$
whenever $\Omega$ is a conormal extension of $\Lambda$ determined by an open convex cone.

\section{Skeleta from toric varieties}
\label{sec:skeleta}

We wish to probe the HMS conjecture \ref{conj:hms3} using the coherent-constructible correspondence of \cite{FLTZ} and \cite{FLTZ2}.  We will apply the theory of those papers to the toric orbifold $\X$ of Section \ref{sec:rpams} and its subvarieties, but let us start by recalling some features of the general case.

Let us first assume that $\cX$ is a toric variety with no stackiness.  Let $\Sigma \subset \cN_\bR$ be the fan associated to $\cX$, where $\cN_\bR = \cN \otimes \bR$ for a suitable lattice $\cN$.  Then, letting $\cM = \Hom(\cN,\bZ)$ and $\cM_\bR = \Hom(\cN_\bR,\bR)$, we (with Fang and Liu) defined in \cite{FLTZ} a conic Lagrangian subset of $T^*(\cM_\bR / \cM)$, as follows:
\[
\Lambda_\Sigma = \bigcup_{\sigma \in \Sigma} (\sigma^\perp + \cM) \times - \sigma
\label{eq:LamSig}
\]

Here $\sigma^\perp \subset \cM_\bR$ denotes the collection of linear covectors on $\cN_\bR$ that vanish on every element of $\sigma \in \Sigma$, and $\sigma^\perp + \cM$
denotes its image in $\cM_\bR/\cM$.  As $\sigma$ is rational, $\sigma^\perp + \cM$ is a subtorus of $\cM_\bR/\cM$, and $(\sigma^\perp + \cM) \times -\sigma$ can be identified with an open subset of its conormal bundle.

When $\cX$ is an orbifold we have the following variant from \cite{FLTZ2}.  Actually we will give a slightly more general treatment here than in \cite{FLTZ2}, and allow $\cX$ to have some ``gerbiness'' or generic isotropy.  From $\cX$ we get a stacky fan in the sense of \cite{BCS}, which we will denote by $\mathbf{\Sigma}$.  In more detail, $\mathbf{\Sigma}$ consists of the following data:
\begin{itemize}
\item A finitely generated abelian group $\cN$.  (In our applications, $\cN$ will be a quotient of the group $\N$ of Section \ref{sec:rpams}.)  
\item A complete simplicial fan $\Sigma \subset \cN_\bR := \cN \otimes_\bZ \bR$
\item A collection $\beta = \{\beta_i\}$, one
$\beta_i$ for each ray (i.e., one-dimensional cone) $\rho_i \in \Sigma$, such that
$\beta_i \in \cN$ maps to a nonzero element $\overline{\beta}_i$
of $\rho_i$ under the natural map $\cN \to \cN_\bR$.
\end{itemize}
To define the conic Lagrangian associated to a stacky fan \cite[Definition 6.3]{FLTZ2}, we replace the ambient torus $\cM_\bR/\cM$ by the Pontrjagin dual of $\cN$, i.e. by the group $G$ of homomorphisms $\cN \to \bR/\bZ$.  If $\cN$ has torsion, then $G$ may be disconnected.  We replace the subtori $\sigma^\perp + \cM_\bZ$ of formula \eqref{eq:LamSig} by closed subgroups (again, possibly disconnected) $G_\sigma \subset G$.  Specifically, for each $\sigma \in \Sigma$, $G_\sigma$ is the group of homomorphisms $\phi:\cN \to \bR/\bZ$ that have $\phi(\beta_i) = 1$ whenever $\beta_i \in \sigma$.  Then we define $\Lambda_{\mathbf \Sigma}$ by
\[
\Lambda_{\mathbf \Sigma} = \bigcup_{\sigma \in \Sigma} G_\sigma \times -\sigma
\label{eq:LambfSig}
\]

The results from \cite{FLTZ,FLTZ2} give us the following:

\begin{theorem}[Coherent-constructible correspondence, or ``CCC'']
\label{thm:fltz}
There is a full embedding
$$\kappa:\Perf(\cX) \hookrightarrow \Sh(G;\Lambda_{\mathbf \Sigma})$$
\end{theorem}

Conjecturally, this embedding is an equivalence---the ``equivariant'' version of this is proven in \cite{FLTZ,FLTZ2}.

\subsection{Understanding the Kashiwara-Schapira sheaf on $\Lambda_{\mathbf \Sigma}$}
\label{sec:understandingKS}

The sheaf structure of $\MSh_{\Lambda_{\mathbf \Sigma}}$ gives us a restriction operation
$$\Sh(G;\Lambda_{\mathbf \Sigma}) \to \MSh_{\Lambda_{\bf \Sigma}}(\Omega)$$
for each open subset $\Omega \subset \Lambda_{\bf \Sigma}$.  For certain $\Omega$ these operations have coherent counterparts in the CCC.  We describe these in this section.

If $\cX$ is a toric variety with fan $\Sigma \subset \cN_\bR$, then the closure of any torus orbit is also a toric variety in a natural way.  The structure of this toric subvariety is visible in the combinatorics of the fan.  Each cone $\sigma \in \Sigma \subset \cN_\bR$ has an associated ``normal fan''
$$\Sigma(\sigma) \subset \cN_\bR/(\bR\cdot \sigma)$$
consisting of the projections of those $\tau \in \Sigma$ with $\tau \supset \sigma$.

These remarks extend to stacky fans.
If $\cN$ is a finitely generated abelian group and ${\mathbf \Sigma} = (\Sigma \subset \cN_\bR,\beta \subset \cN)$ is a stacky fan, then for each $\sigma$ we may form the quotient $\cN(\sigma) = \cN/(\bZ\cdot \{\beta_i : \overline{\beta}_i \in \sigma\})$.  Then the fan $\Sigma(\sigma)$ lives naturally in $\cN(\sigma)_\bR$, and the subcollection
$\beta_\sigma = \{\beta_j : \overline{\beta}_j
\in \tau\setminus \sigma \hbox{ for some } \tau\supset \sigma\}$
of $\beta_j$ with
$\overline{\beta}_j$ in a cone containing
$\sigma$ but not in $\sigma$ itself,
forms a stacky structure on this fan.  We denote the resulting stacky fan by ${\mathbf{\Sigma}}(\sigma)$.  There is a natural inclusion of toric stacks $\cX_{{\mathbf \Sigma}(\sigma)} \hookrightarrow \cX_{\mathbf{\Sigma}}$.

The skeleton $\Lambda_{{\mathbf \Sigma}(\sigma)}$ is a conic Lagrangian in the cotangent bundle of the subgroup $G_\sigma \subset G$.  We define an open subset of $\Lambda_{\mathbf \Sigma}(\sigma) \subset \Lambda_{\mathbf \Sigma}$ as follows:
$$\Lambda_{\mathbf \Sigma}(\sigma) = \bigcup_{\tau \supset \sigma} G_\tau \times (-\tau^\circ)$$
This is just the conormal extension of $\Lambda_{\Sigma(\sigma)}$ by the trivial bundle of cones determined by $\sigma^\circ$ (see Section \ref{sec:conext}).

\begin{theorem}
\label{thm:understandingKS}
We have a commutative square
$$
\xymatrix{
\Perf(\cX_{\mathbf \Sigma}) \ar[d]_{\kappa} \ar[rr] & & \Perf(\cX_{{\mathbf \Sigma}(\sigma)} \ar[d]^{\kappa}) \\
\MSh_{\Lambda_{\mathbf\Sigma}}(\Lambda_{\mathbf\Sigma}) \ar[r]_{\mathit{res}\quad} & \MSh_{\Lambda_{\mathbf\Sigma}}(\Lambda_{\mathbf\Sigma}(\sigma)) \ar[r]_{\cong\,} & \MSh_{\Lambda_{{\mathbf \Sigma}(\sigma)}}(\Lambda_{{\mathbf\Sigma}(\sigma)})
}
$$
The top row is given by restriction of vector bundles, the columns are given by the CCC, and the bottom row is given by restriction in the Kashiwara-Schapira sheaf and conormal extension.
\end{theorem}

\section{Skeleta from toric hypersurfaces}
\label{sec:stackyfan}

Now we return to the Batyrev-Borisov hypersurfaces.  Let $\triangle,\triangle^\vee$, $\X,\Xvee$ and $\Z^\infty,\Zvee_\infty$ be as in Section \ref{sec:rpams}.  We assume that $\triangle$ is vertex-simplicial, and $\triangle^\vee$ is facet-simplicial, so that in particular we regard $\X$ as a smooth Deligne-Mumford stack.  In more detail, $\X$ is determined by a stacky fan $\bSi = (\Sigma,\{\beta_i\})$ that is closely related to $\triangle^\vee \subset \Mvee = N$.  In the notation of Section \ref{sec:skeleta}, we have $\cN = N$, a free abelian group.  The cones of  $\Sigma$ are indexed by the faces of $\triangle^\vee$, by setting $\sigma_F := \bR_{\geq 0} \cdot F \subset N_\bR$, and the $\beta_i$ are precisely the vertices of $\triangle^\vee$.

Now we attach to $\bSi$ the conic Lagrangian $\Lambda = \Lambda_{\bSi}$ as in Section \ref{sec:skeleta}.  As $\M_\bR = \Nvee_\bR$ and $\M = \Nvee,$
after choosing an inner product we may identify
$T^*(\M_\bR/M)$ with the complex algebraic torus $\Tvee$.  The CCC and microlocalization suggest that $\Lambda \subset \Tvee$ is the skeleton of a Weinstein
structure\footnote{A Weinstein manifold is a complete, exact symplectic manifold
whose Liouville vector
field is gradient-like for an exhausting Morse function.  The skeleton is the union of stable manifolds of Liouville flow.} on $\Tvee$ that is suited for studying the mirror of $X$.  We propose that the contact boundary of $\Lambda$, i.e. $\Lambda^\infty \subset \M_\bR/\M \times (\N_\bR - 0)/(\bR_{>0})$, is suited for studying the mirror $\Zvee_\infty$ of $\Z^\infty$.  In particular $\Lambda^\infty$ should be homeomorphic to a skeleton of $\Zvee_\infty$.

\begin{conjecture}
\label{conj:skeletonconjecture}
Suppose $\triangle$ is vertex simplicial and $\triangle^\vee$ is facet simplicial.  There is a
Weinstein structure on $\Zvee_\infty$ whose Lagrangian skeleton is homeomorphic to $\Lambda^\infty$.  The sheaf of Fukaya categories on this skeleton is equivalent to the Kashiwara-Schapira sheaf on $\Lambda^\infty$.
\end{conjecture}

In view of mirror symmetry, the following result is strong evidence for the Conjecture.

\begin{theorem}
\label{thm:52}
The category of global sections of the Kashiwara-Schapira sheaf on $\Lambda^\infty$ contains $\Perf(\Z^\infty)$ as a full subcategory.
\end{theorem}

A proof will appear in forthcoming work with Sibilla \cite{STZ-inprogress}.  As discussed in Section \ref{sec:understandingKS}, the normal stacky fans $\bSi(\sigma)$ determine an open cover $\{\Lambda_{\bSi}(\sigma)\}$ of $\Lambda_\bSi$.  If we omit the zero cone, we get a closely related open cover $\{\Lambda_{\bSi}(\sigma)/\bR_{>0}\}$ cover $\Lambda^\infty$ which can be used to compute the global sections of the Kashiwara-Schapira sheaf.  Theorems \ref{thm:understandingKS} yields the full embedding Theorem \ref{thm:52} by applying the CCC to each component of $\Z^\infty$ and each open chart of $\Lambda^\infty$.  This full embedding is an equivalence as long as the conjecture alluded to after Theorem \ref{thm:fltz} holds.

\subsection*{Concrete Description of $\Lambda^\infty$}
In the notation of Section \ref{sec:skeleta}, we begin with $\cN = N,$ a free abelian group.
Each face $F$ of $\triangle^\vee$ labels a cone $\sigma_F := \bR_{\geq 0}\cdot F\subset N_\bR.$
The stacky fan $\bSi(\sigma_F)$ has associated abelian group $N_F := \cN(\sigma_F)
:= N/(\bZ\cdot V_F),$ not necessasrily free; here $V_F$ is the set of vertices of $F.$
Then $M_F = V_F^\perp.$
Cones of $\bSi(\sigma_F)$ are labeled by faces $F'\supset F$
and defined by $\sigma^F_{F'} = \bR_{\geq 0}\cdot V_{F'},$ or rather the image of this set in
$\cN(\sigma_F)_\bR = N_\bR/(\bR\cdot V_F).$
The set of lifts of ray vectors is $\beta_F = \bigcup_{F'\supset F}V_{F'}\setminus V_F.$
The associated abelian group is $G_{F'}:= G_{\sigma_{F'}} =
\{m\in M_\bR : m(V_{F'})\subset \bZ\}/M\subset M_\bR/M.$
Note that this group depends only on $F'$, not the pair $(F',F).$
We have a short exact sequence of ableian groups
$1\rightarrow G^\circ_{{F'}}\rightarrow
G_{{F'}}\rightarrow P_{{F'}}\rightarrow 1,$
where $G^\circ_{F'}:= M_{F',\bR}/M = V_{F'}^\perp/M$ is the neutral component
of $G_{F'},$ a torus, and $P_{F'}=\pi_0(G_{F'})$ is a finite
quotient.  Then $\Lambda' := \bigcup_F G_F\times F^\circ$ is the union of conormal
extensions of the skeleta for the toric components which comprise $Z^\infty.$
We define $\Lambda^\infty = \Lambda'/\bR_{>0},$ which has an equivalent microlocal
sheaf.

\begin{figure}[ht]
\begin{center}
\includegraphics[scale=0.65]{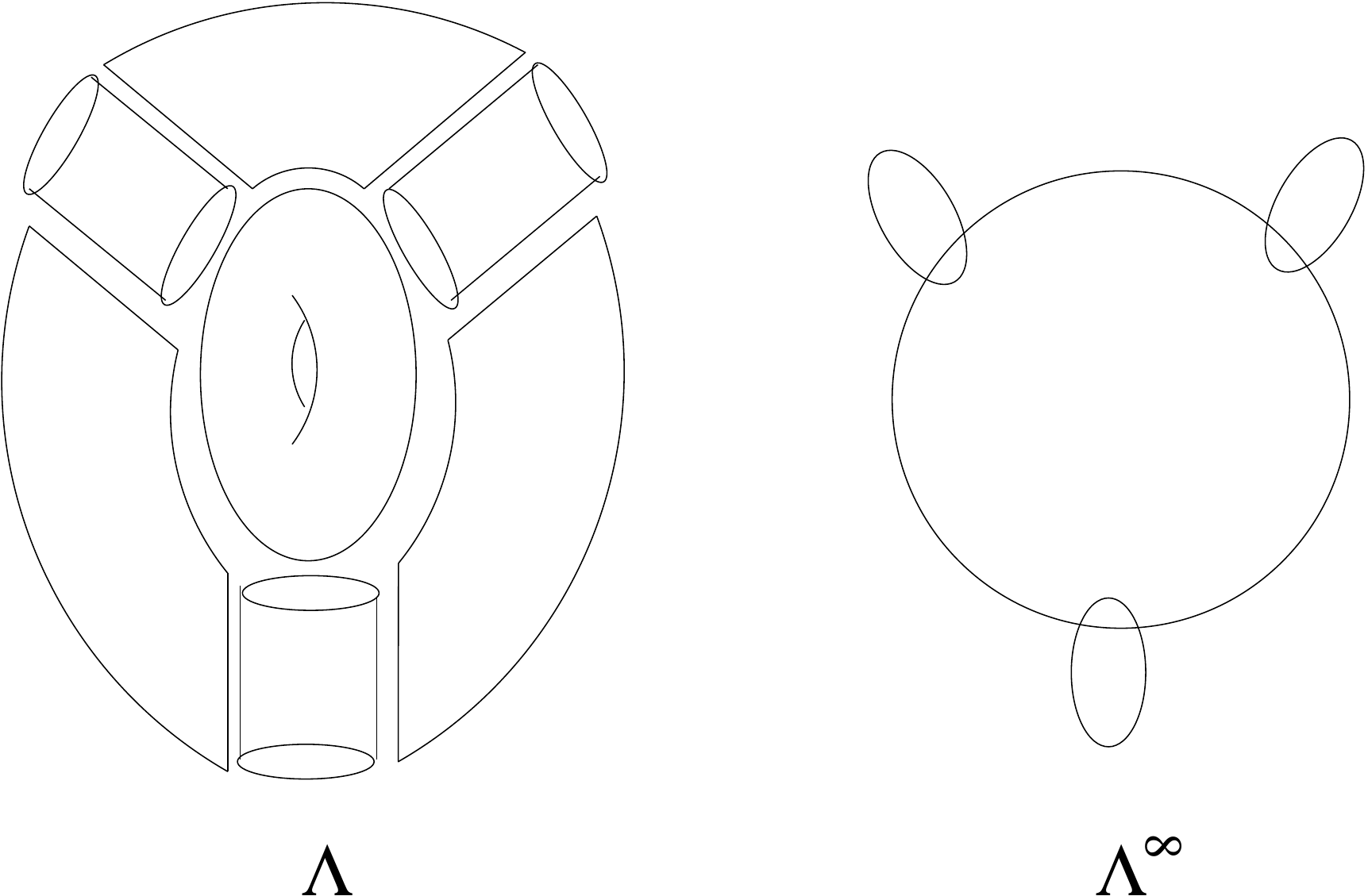}
\end{center}
\caption{On the left is a rough picture of how the conic Lagrangian $\Lambda$ attached to $\bP^2$ should be constructed from constituents $G_\sigma \times - \sigma$.  
On the right is the Legendrian $\Lambda^\infty$ that we expect to be a skeleton of $\Zvee_\infty = \{(x,y) \in \bC^* \times \bC^* \mid a x + b y + c(xy)^{-1} + d = 0\}$.}
\end{figure}

\section{The homology of $\Lambda^\infty$}
\label{sec:cc}

A consequence of Conjecture \ref{conj:skeletonconjecture} is that
we have  a
homotopy equivalence between $\Zvee_\infty$ and $\Lambda^\infty$. 
A weaker conjecture then states that these spaces have the same rational homology. 
We will discuss a computer program that
makes many tests of this weaker conjecture in Section \ref{subsec:compprog}.

The space $\Lambda^\infty$ projects naturally onto $(\N_\bR - 0)/\bR_{>0}$, which we can identify with $\partial \triangle^\vee$.  This induces a Leray spectral sequence abutting to the homology
of $\Lambda^\infty.$
The projection is trivial over the interior of each face of $\partial \triangle^\vee$, so that
$$\Lambda^\infty = \bigcup_F G_F \times F^\circ$$
where $G_F := G_{\sigma_F}$ are the disconnected tori defined in Section \ref{sec:skeleta} and the union is over proper faces of $\triangle^\vee$---i.e. $G_F$ is the group of maps from $\N$ to $\bR/\bZ$ that carries each $\beta_i \in F$ to $1$.  

The natural inclusion $\phi^F_{F'}:G_{F'} \hookrightarrow G_F$ for $F \subset F'$ describes the attaching of $G_F \times F^\circ$ to $G_{F'} \times F'^\circ$.  The $E^1$
page of the Leray spectral sequence is the bi-graded chain complex
$$
0\leftarrow
\bigoplus_{\dim\,F = 0} H_*(G_F) \leftarrow \bigoplus_{\dim\,F = 1} H_*(G_F) \leftarrow \bigoplus_{\dim\,F = 2} H_*(G_F)\leftarrow\cdot\cdot\cdot
$$
whose differentials are alternating sums of the maps induced by $\phi^F_{F'}$.  We conjecture that
with rational coefficients the spectral sequence degenerates, i.e. $E^2 = E^\infty.$

Explicitly, we define a chain complex $(C_*,\partial)$ as follows.  For each face $F$, let $P_F = \pi_0(G_F)$ and let $G^\circ_F \subset G_F$ be the neutral component of $G_F$, so that we have a short exact sequence
$$1 \to G_F^\circ \to G_F \to P_F \to 1,$$
with $G_F^\circ = M_{F,\bR}/M_F$ (or $M_F = H_1(G^\circ_F)$.  We have canonically $H_*(G_F) = (\bigwedge^* M_F)^{\oplus |P_F|}$.  Pick an orientation of each face of $\triangle^\vee$, so that for each pair of faces $F',F$ with $F$ of codimension one in $F'$ we get a sign $\epsilon_{F,F'} = \pm 1$.  Then for each $\gamma \in (\bigwedge^r M_F)^{\oplus |P_F|}$ we set
$$\partial \gamma = \sum_{F'} \epsilon_{F,F'} (\phi^F_{F'})_* (\gamma)$$
If $\partial \gamma = 0$ then $\gamma$ contributes to $H_{\dim(F) + r}(\Lambda^\infty)$.

\subsection{Example}

Consider the triangle $\triangle^\vee = \mathit{conv}\left(\{(2,-1),(-1,2),(-1,-1)\}\right)
\subset \bR^2 = M^\vee \otimes_\bZ \bR,$ where $M^\vee = \bZ^2.$ 
Let $F = \{(-1,-1)\}$ so $\sigma_F = \bR_{\geq 0}\cdot \{(-1,-1)\}$ and $\cN(\sigma_F) = 
\bZ^2/\{(-1,-1)\} \cong \bZ,$ with $(a,b)\sim (a-b)\in \bZ.$
Then if $F'\supset F$ is the edge containing $(2,-1)$ and $(-1,-1),$
$\sigma^F_{F'}$ is the $\bR_{\geq 0}$ span of $(3,0)\sim 3\in \bZ.$ 
Let $F''\supset F$ be the remaining edge.  Then
$\sigma^F_{F''}$ 
is the span of $-3\in \bZ.$  Further, $\beta_F = \{3,-3\}$
and this gives the stacky fan structure for the toric stack $\bP^1(3,3),$
a $\bP^1$ with $\bZ/3\bZ$ orbifold points at the two poles.
$G_{F} = \{(a,b) : -a-b \in \bZ\}/\bZ^2 \cong S^1,$ so $H_*(G_F)\cong \bZ\oplus \bZ[-1].$
A similar story holds at each vertex.  For the edge $F'$ we have $G_{F'}=
\{(a,b): -a-b\in \bZ, 2a-b\in \bZ\}/\bZ^2\cong \bZ/3\bZ.$
We can label the three points $(k/3,1-k/3),$ and note that under
$G_{F'}\hookrightarrow G_F$ they are distributed evenly on $S^1,$ and this explains
how to attach chains to make $\Lambda^\infty.$  As chains in $\Lambda^\infty,$
these three points contribute in dimension $0+1=1,$ since they lie over edges.
A similar story holds at each edge.

The space of chains looks like $\bQ^3$ in degree zero (coming from the vertices)
and $\bQ^3 \oplus \bQ^9$ in degree one (the first factor coming from the vertices,
the second from the edges).  The only nonzero part of the differential is its restriction
to the edges: $\bQ^9\rightarrow \bQ^3,$ and its form comes purely from the combinatorics
of the polygon $\triangle^\vee,$
repeated three times because of the $\bZ/3\bZ$'s.
The $3\times 9$ matrix representing this map is
 $$\begin{pmatrix}
-1 & -1 & -1 & -1 & -1 & -1 & 0 & 0 & 0 \\
1 & 1 & 1 & 0 & 0 & 0 & -1 & -1 & -1 \\
0 & 0 & 0 & 1 & 1 & 1 & 1 & 1 & 1 \\
\end{pmatrix}.$$
We conclude that the image is two-dimensional, and hence $\dim H_0(C_*)
 = 3-2 = 1,$
i.e. $H_0(C_*) = H_0(\Lambda^\infty) \cong \bQ$
while $H_1(C_*) = H_1(\Lambda^\infty) \cong \bQ^3\oplus \bQ^7.$

The polygon $\triangle^\vee$ defines the line bundle $\mathcal O_{\bP^2}(3),$ and the
affine cubic is an elliptic curve minus nine points at infinity (three for each toric divisor).
It is connected with Euler characteristic $0-9,$ hence $\dim H_0(Z^\vee_\infty) = 1$ and
$\dim H_1(Z^\vee_\infty) = 10,$
which agrees with our computation.

\subsection{Computer program}
\label{subsec:compprog}
It is straightforward to automate the process of calculating the homology of $(C_*,\partial)$
for any polytope, including for example
all facet-simplicial reflexive polytopes in dimension three.
We have done this and posted the Sage worksheet in a public folder
at {\tt http://dl.dropbox.com/u/24939613/skeleton-homology.sws}
(entering the URL in a browser will trigger a download of the file).
Here we explain the algorithm.
The program code and the example cell are also generously annotated.

The program runs by calculating the ranks and kernels of matrices encoding the
maps $(\varphi^F_{F'})_*$ on the $\bZ$-modules $H_*(G_F,\bZ)$. 
Since $H_*(G_F) \cong (\bigwedge^*M_F)^{\oplus |P_F|},$
after we choose a $\bZ$-basis $e_i,$ $i \in \{1,...,s\},$ for $M_F$ this determines a
basis $e_I^k$ for $\bigwedge^*M_F$, one basis element for each subset $I\subset \{1,...,s\}$ and element
$k\in P_F,$ or $2^s|P_F|$ in all.

First choose an identification $N\cong \bZ^n$ and suppose $\triangle^\vee
= \mathit{conv}\{v_1,...,v_p\},$
with $v_i = (v_{i,1},...,v_{i,n}).$  A face $F$ is indexed by a subset $J = \{j_1 < j_2 < \cdots\} \subset \{1,\ldots, p\}$ indexing its vertices.  Let $A_J$ be the matrix with $(A_J)_{a,b} = v_{j_a,b}$.  
By changing bases for the lattice $\bZ\cdot V_F$ and for $N$,
we can perform row \emph{and} column operations to put
$A_J$ in Smith normal form (computed by {\tt smith\_form()} in Sage)
$$
S_J = U_J A_J V_J,
$$
where $S_J$ is diagonal with entries $s_{11},\ldots, s_{rr}$
These matrices fix an identification of $G_F$ with
$(\bR/\bZ)^s \times \bZ/s_{11} \times \cdots \times \bZ/s_{rr}$.

If $F$ is a face corresponding to $J \subset \{1,\ldots,p\}$ and $F' \supset F$
is a face of one dimension higher corresponding to
$K \supset J$, the matrix $V_{J}^{-1} V_K$ encodes the information of the
map $\phi^F_{F'}$. 

The southeastern $(s+1)\times s$ block represents
a linear map from $M_{F'}$ to $M_F.$  The linear
maps from $\Lambda^*M_{F'}$ to $\Lambda^*M_F$ are described by all the various minors of
this block.  The northwestern $(r-1)\times r$ block of $V_J^{-1}V_K$ encodes the map
of discrete groups:  to find out where the element $(k_1,...,k_r)\in P_{F'}$ goes, act on it
by the matrix and read off the answer as an element of the group $P_F.$  (Note that in
restricting to these blocks, we are ignoring some
information, \emph{irrelevant to homology}, in the matrix $V_J^{-1}V_K$
which describes where the discrete group $P_{F'}$ (in the induced splitting) maps
into the group $G_{F}.$

With these data in hand, the only other information needed to implement the homology
calculation is the differential of the polytope, and that is accomplished through simple
combinatorics.

\subsection{Comparisons}
Danilov and Khovanski\v{\i} \cite{DK} have computed the ranks of the homology groups
(or rather the dual compactly supported cohmology) of affine hypersurfaces in
toric varieties.  We have compared dimensions
and verified that our chain complex computes the right betti numbers, as
conjectured, for all 194 facet-simplicial three-dimensional
reflexive polytopes.

\end{document}